\documentclass[12pt]{amsart}
\usepackage{amsfonts}
\usepackage{latexsym}
\usepackage{amssymb}
\usepackage{amsmath}

\renewcommand{\H}{\mathbb{H}}
\newcommand{\B}{\mathbb{B}}

\newcommand{\N}{\mathbb{N}}

\newcommand{\cC}{\mathcal{C}}
\newcommand{\cD}{\mathcal{D}}

\newcommand{\cF}{\mathcal{F}}

\newcommand{\cH}{\mathcal{H}}

\newcommand{\cN}{\mathcal{N}}

\newcommand{\cP}{\mathcal{P}}

\newcommand{\cS}{\mathcal{S}}

\newcommand{\cs}{\mathfrak s}

\newcommand{\ep}{\varepsilon}

\newcommand{\sm}{\setminus}

\newcommand{\Lb}{{\Big\lbrace}}
\newcommand{\Rb}{{\Big\rbrace}}
\newcommand{\lb}{{\big\lbrace}}
\newcommand{\rb}{{\big\rbrace}}
\newcommand{\ls}{\mbox{\large $($}}
\newcommand{\rs}{\mbox{\large $)$}}

\newcommand{\diams}{\mbox{\rm diam}\;\!}
\newcommand{\diam}{\mbox{\rm diam}}

\newcommand{\bcup}{\bigcup}
\newcommand{\mcup}{\mbox{\small $\bigcup$}}

\newcommand{\res}{\mbox{\LARGE{$\llcorner$}}}

\newcommand{\eR}{{\overline {\mathbb R}}}

\newtheorem{The}{Theorem}

\newtheorem{Def}{Definition}
\newtheorem{Rem}{Remark}

\begin{document}

\title
[On a measure-theoretic area formula]
{{\bf On a measure-theoretic area formula}}
\author{Valentino Magnani}
\address{Valentino Magnani, Dipartimento di Matematica \\
Largo Pontecorvo 5 \\ I-56127, Pisa}
\email{magnani@dm.unipi.it}
\date{\today}
%
%
%
%

\thanks{The author acknowledges the support of the European Project ERC AdG *GeMeThNES*.}

\begin{abstract}
We show how classical differentiation theorems for measures can be turned into an integral representation of a Borel measure with respect to a fixed Carathéodory measure. We focus our attention on the cases where this measure is both the Hausdorff measure and the spherical Hausdorff measure, giving the corresponding measure-theoretic area formula. Our point consists in using certain covering derivatives as ``generalized densities''. Some consequences for the sub-Riemannian Heisenberg group are also pointed out.
\end{abstract}

\maketitle

It is well known that computing the Hausdorff measure of a set by an integral formula
is usually related to rectifiability properties, namely, the set must be close to a linear subspace at any small scale. The classical area formula exactly provides this relationship by an integral representation of the Hausdorff measure. Whenever a rectifiable set is thought of as a countable union of Lipschitz images of subsets in a Euclidean space, the area formula holds in metric spaces, \cite{Kir94}.

In the last decade, the development of Geometric Measure Theory in a non-Euclidean framework raised new 
theoretical questions on rectifiability and area-type formulae. The main problem in this setting stems from 
the gap between the Hausdorff dimension of the target and that of the source space of the parametrization.
In fact, in general this dimension might be strictly greater than the topological dimension of the set.
As a result, the parametrization from a subset of the Euclidean space cannot be Lipschitz continuous
with respect to the Euclidean distance of the source space.

To mention an instance of this difficulty, the above mentioned area formula for a large class 
of Heisenberg group valued Lipschitz mappings does not work, \cite{AmbKir00}.
For this reason, theorems on differentiation of measures constitute an important tool to overcome this problem. 
In this connection, the present note shows how the Federer's Theorems of 2.10.17 and 2.10.18 in \cite{Federer69}
are able to disclose a purely metric area formula. The surprising aspect of this formula is that
an ``upper covering limit'' actually can be seen as a {\em generalized density} of a fixed Borel measure.

To define these densities, we first introduce covering relations: if $X$ is any set, a {\em covering relation} is a subset $\cC$ of $\{(x,S): x\in S\in\cP(X)\}$. 
In the sequel, the set $X$ will be always assumed to be equipped with a distance.
Defining for $A\subset X$ the corresponding class  $\cC(A)=\{S: x\in A,\,(x,S)\in\cC\}$, 
we say that $\cC$ is {\em fine at $x$}, if for every $\ep>0$ there exists $S\in\cC(\{x\})$ such that $\diams S<\ep$.
According to 2.8.16 of \cite{Federer69}, the notion of covering relation yields the following notion
of ``covering limit''.
%
%
%
%
%
\begin{Def}[Covering limit]\label{def:coveringlimits} \rm
If $\cC$ denotes a covering relation and $f:\cD\to\eR$, $\cC(\{x\})\subset \cD\subset\cC(X)$ and $\cC$ is fine at $x\in X$, then we define the {\em covering limits}
\begin{eqnarray}
&&(\cC)\limsup_{S\to x}f=\inf_{\ep>0}\sup\{f(S): S\in\cC(\{x\}), \diams S<\ep\}\,,\\ 
&&(\cC)\liminf_{S\to x}f=\sup_{\ep>0}\inf\{f(S): S\in\cC(\{x\}), \diams S<\ep\}\,.
\end{eqnarray}
\end{Def}
The covering relations made by closed balls clearly play an important role
in the study of the area formula for the spherical Hausdorff measure.
\begin{Def}\rm
The closed ball and the open ball of center $x\in X$ and radius $r>0$ are denoted by 
\[ 
\B(x,r)=\{y\in X: d(x,y)\le r\}\quad\mbox{and}\quad B(x,r)=\{y\in X: d(x,y)< r\}\,,
\]
respectively. We denote by $\cF_b$ the family of all closed balls in $X$.
\end{Def}
The next definition introduces the Carath\'eodory construction, see 2.10.1 of \cite{Federer69}.
%
%
%
%
\begin{Def}\rm 
Let $\cS\subset\cP(X)$ and let $\zeta:\cS\to[0,+\infty]$ represent the {\em size function}.
If $\delta>0$ and $R\subset X$, then we define
\begin{eqnarray*}
&&\zeta_\delta(R)=\inf \Big\lbrace\sum_{j=0}^\infty \zeta(E_j):
E_j\in\cF,\ \diam(E_j)\le\delta\ \mbox{for all}\ j\in\N, R\subset \bcup_{j\in\N} E_j \Big\rbrace\,.
\end{eqnarray*}
The {\em $\zeta$-approximating measure} is defined and denoted by $\psi_\zeta=\sup_{\delta>0}\zeta_\delta$.
Denoting by $\cF$ the family of closed sets of $X$, for $\alpha,\,c_\alpha>0$, we define
$\zeta_\alpha:\cF\to[0,+\infty]$ by
\begin{equation*}
\zeta_\alpha(S)= c_\alpha\, \diam(S)^\alpha\,.
\end{equation*}
Then the {\em $\alpha$-dimensional Hausdorff measure} is $\cH^\alpha=\psi_{\zeta_\alpha}$. 
If $\zeta_{b,\alpha}$ is the restriction of $\zeta_\alpha$ to $\cF_b$, 
then $\cS^\alpha=\psi_{\zeta_{b,\alpha}}$ is the
{\em $\alpha$-dimensional spherical Hausdorff measure}.
\end{Def}
These special limits of Definition~\ref{def:coveringlimits} naturally arise 
in the differentiation theorems for measures
and allow us to introduce a special ``density'' associated with a measure.
%
%
%
%
%
\begin{Def}[Federer density]\rm 
Let $\mu$ be a measure over $X$, let $\cS\subset\cP(X)$ and let $\zeta:\cS\to[0,+\infty]$. 
Then we set 
\[ 
\cS_{\mu,\zeta}=\cS\sm \Big\{S\in\cS\,:\; \zeta(S)=\mu(S)=0\; \mbox{\rm or}\; 
\mu(S)=\zeta(S)=+\infty\, \Big\}\,,
\]
along with the covering relation $\cC_{\mu,\zeta}=\{(x,S): x\in S\in\cS_{\mu,\zeta}\}$.
We choose $x\in X$ and assume that $\cC_{\mu,\zeta}$ is fine at $x$. We define the quotient function
\begin{eqnarray*}
Q_{\mu,\zeta}:\cS_{\mu,\zeta}\to[0,+\infty],\quad
Q_{\mu,\zeta}(S)=\left\{\begin{array}{ll} +\infty & \mbox{if $\zeta(S)=0$} \\ 
\mu(S)/\zeta(S) & \mbox{if $0<\zeta(S)<+\infty$} \\
0 & \mbox{if $\zeta(S)=+\infty$}
\end{array}\right.\,.
\end{eqnarray*}
Then we are in the position to define the {\em Federer density}, or 
{\em upper $\zeta$-density of $\mu$ at $x\in X$}, as follows
\begin{equation}\label{eq:thetazeta}
F^\zeta(\mu,x)=(\cC_{\mu,\zeta})\limsup_{S\to x}Q_{\mu,\zeta}(S)\,.
\end{equation}
\end{Def}
According to the following definition, we will use special notation 
when we consider Federer densities with respect to 
$\zeta_\alpha$ and $\zeta_{b,\alpha}$, respectively.
\begin{Def}\rm 
If $\mu$ is a measure over $X$ and $\cC_{\mu,\zeta_{b,\alpha}}$ is fine at $x\in X$, then we set
$\theta^\alpha(\mu,x)=F^{\zeta_{b,\alpha}}(\mu,x)$. If $\cC_{\mu,\zeta_\alpha}$ is fine at $x$,
then we set $\cs^\alpha(\mu,x)=F^{\zeta_\alpha}(\mu,x)$.
\end{Def}
\begin{Rem}\rm
If $x\in X$ and there exists an infinitesimal sequence $(r_i)$ of positive radii such that 
$\B(x,r_i)$ may have vanishing diameter and in this case $\mu\ls\B(x,r_i)\rs>0$,
then it is easy to realize that both $\cC_{\mu,\zeta_{b,\alpha}}$ and $\cC_{\mu,\zeta_\alpha}$ are fine at $x$. 
In particular, the last conditions are always satisfied whenever all balls $\B(x,r_i)$ have positive diameter.
\end{Rem}
By the previous definitions, we state a revised version of Theorem~2.10.17(2) of \cite{Federer69}.
%
%
%
%
%
%
%
%
%
%
%
%
%
%
%
\begin{The}\label{the:Idiffthegauge'}
Let $\cS\subset\cP(X)$ and let $\zeta:\cS\to[0,+\infty]$ be a size function.
If $\mu$ is a regular measure over $X$, $A\subset X$, $t>0$, $\cS_{\mu,\zeta}$ covers $A$ finely
and for all $x\in A$ we have $F^\zeta(\mu,x)<t$, 
then $\mu(E)\le t\,\psi_\zeta(E)$ for every $E\subset A$.
\end{The}

%
%
%
%
%
%
%
%
%
%
%
%
%
%
Analogously, the next theorem is a revised version of Theorem~2.10.18(1) in \cite{Federer69}.
\begin{The}\label{the:IIdiffthegauge}
Let $\mu$ be a measure over $X$, let $\cS$ be a family of closed and $\mu$-measurable sets,
let $\zeta:\cS\to[0,+\infty)$, let $B\subset X$ and assume that $\cS_{\mu,\zeta}$ covers
$B$ finely. If there exist $c,\eta>0$ such that for each $S\in\cS$ there exists $\tilde S\in \cS$ with the properties 
\begin{equation}\label{eq:hatStildeS}
\hat S\subset\tilde S,\quad \diams\tilde S\le c\,\diams S\quad\mbox{and}\quad \zeta(\tilde S)\le\eta\,\zeta(S),
\end{equation}
where $\hat S=\mcup\,\{ T\in\cS: T\cap S\neq\emptyset, \diams T\le2\,\diams S\}$, 
$V\subset X$ is an open set containing $B$ and for every $x\in B$ we have 
$F^\zeta(\mu,x)>t$, then $\mu(V)\ge t\,\psi_\zeta(B)$.
\end{The}

These theorems provide both upper and lower estimates for a large class of measures, 
starting from upper and lower estimates of the Federer density.
A slight restriction of the assumptions in the previous theorems joined with some standard arguments
leads us to a new metric area-type formula, where the integration of $F^\zeta(\mu,\cdot)$ 
recovers the original measure. This is precisely our first result.
%
%
%
%
%
%
%
%
\begin{The}[Measure-theoretic area-type formula]\label{the:meastheoarea}
Let $\mu$ be a Borel regular measure over $X$ such that there exists a countable open covering of $X$,
whose elements have $\mu$ finite measure. Let $\cS$ be a family of closed sets, 
let $\zeta:\cS\to[0,+\infty)$ and assume that for some constants $c,\eta>0$ and for every $S\in\cS$ 
there exists $\tilde S\in \cS$ such that  
\begin{equation}\label{eq:hatSdoubling}
\hat S\subset\tilde S,\quad \diams\tilde S\le c\,\diams S\quad\mbox{and}\quad \zeta(\tilde S)\le\eta\,\zeta(S),
\end{equation}
where $\hat S=\mcup\,\{ T\in\cS: T\cap S\neq\emptyset, \diams T\le2\,\diams S\}$.
If $A\subset X$ is Borel, $\cS_{\mu,\zeta}$ covers $A$ finely and $F^\zeta(\mu,\cdot)$ is a Borel function 
on $A$ with 
\begin{equation}\label{eq:thetalessinf}
\psi_\zeta\ls\lbrace x\in A: \theta^\zeta(\mu,x)=0\rbrace\rs<+\infty\quad\mbox{and}\quad
\mu\ls\lbrace x\in A: \theta^\zeta(\mu,x)=+\infty\rbrace\rs=0\,,
\end{equation}
then for every Borel set $B\subset A$, we have
\begin{equation}\label{eq:mestheoarea}
 \mu(B)=\int_B F^\zeta(\mu,x)\,d\psi_\zeta(x)\,.
\end{equation}
\end{The}
The second condition of \eqref{eq:thetalessinf} precisely corresponds
to the absolute continuity of $\mu\res A$ with respect to $\psi_\zeta\res A$.
This measure-theoretic area formula may remind of a precise differentiation theorem,
where indeed the third condition of \eqref{eq:hatSdoubling} represents a kind of ``doubling condition''
for the size function $\zeta$. In fact, the doubling condition for a measure allows for obtaining 
a similar formula, where the density is computed by taking the limit of the ratio between the measures
of closed balls with the same center and radius, see for instance Theorems 2.9.8 and 2.8.17 of \cite{Federer69}.

On one side, the Federer density $F^\zeta(\mu,x)$ may be hard to compute, depending on the space $X$.
On the other side, formula \eqref{eq:mestheoarea} neither requires special geometric properties for $X$, 
as those for instance of the Besicovitch covering theorem (see the general condition 2.8.9 of \cite{Federer69}),
nor an ``infinitesimal'' doubling condition for $\psi_\zeta\res A$, as in 2.8.17 of \cite{Federer69}.
Moreover, there are no constraints that prevent $X$ from being infinite dimensional.

The absence of specific geometric conditions on $X$ is important especially in relation with applications 
of Theorem~\ref{the:meastheoarea} to sub-Riemannain Geometry, in particular for the class of
the so-called Carnot groups, where the classical Besicovitch covering theorem may not hold, see \cite{KorRei95}, \cite{SawWhe92}. In these groups, we have no general theorem to ``differentiate'' an arbitrary Radon measure,
therefore new differentiation theorems are important.

We provide two direct consequences of Theorem~\ref{the:meastheoarea}, that correspond to the
cases where $\psi_\zeta$ is the Hausdorff measure and the spherical Hausdorff measure, respectively. 
%
%
%
%
%
%
%
%
\begin{The}[Area formula with respect to the Hausdorff measure]\label{the:hausdarea}
Let $\mu$ be a Borel regular measure over $X$ such that there exists a countable open covering of $X$,
whose elements have $\mu$ finite measure.
If $A\subset X$ is Borel and $\cS_{\mu,\zeta_\alpha}$ covers $A$ finely, then
$\cs^\alpha(\mu,\cdot)$ is Borel. Moreover, if $\cH^\alpha(A)<+\infty$ and 
$\mu\res A$ is absolutely continuous with respect to $\cH^\alpha\res A$, then for every Borel set $B\subset A$, we have
\[
 \mu(B)=\int_B \cs^\alpha(\mu,x)\,d\cH^\alpha(x)\,.
\]
\end{The}
This theorem essentially assigns a formula to the density of $\mu$ with respect to $\cH^\alpha$.
Let us recall the formula for this density
\[ 
\cs^\alpha(\mu,x)=\inf_{\ep>0}\Lb\sup\lb Q_{\mu,\zeta_\alpha}(S): x\in S\in\cS_{\mu,\zeta_\alpha}, \diams S<\ep\rb\Rb\,.
\]
Under the still general assumption that all open balls have positive diameter, we have
$Q_{\mu,\zeta_\alpha}(S)=\mu(S)/\zeta_\alpha(S)$. More manageable formulae for 
$\cs^\alpha(\mu,\cdot)$ turn out to be very hard to be found and this difficulty is related
to the geometric properties of the single metric space.
On the other hand, if we restrict our attention to the spherical Hausdorff measure, then the corresponding density $\theta^\alpha(\mu,\cdot)$ can be explicitly computed in several contexts, where it can be also
given a precise geometric interpretation. 

In this case, we will also assume a rather weak condition on the diameters of open balls. 
Precisely, we say that a metric space $X$ is {\em diametrically regular} if for all $x\in X$ and $R>0$ 
there exists $\delta_{x,R}>0$ such that $(0,\delta_{x,R})\ni t\to\diam\ls B(y,t)\rs$ is continuous 
for every $y\in \B(x,R)$.
We are now in the position to state the measure-theoretic area-type formula for the spherical 
Hausdorff measure.

%
%
%
%
%
%
%
\begin{The}[Spherical area formula]\label{the:metricspherical}
Let $X$ be a diametrically regular metric space, let $\alpha>0$ and let $\mu$ be a Borel regular measure over $X$ 
such that there exists a countable open covering of $X$ whose elements have $\mu$ finite measure. 
If $B\subset A\subset X$ are Borel sets and $\cS_{\mu,\zeta_{b,\alpha}}$ covers $A$ finely,
then $\theta^\alpha(\mu,\cdot)$ is Borel on $A$. In addition, if $\cS^\alpha(A)<+\infty$
and $\mu\res A$ is absolutely continuous with respect to $\cS^\alpha\res A$, then we have
\begin{equation}\label{eq:spharea}
 \mu(B)=\int_B \theta^\alpha(\mu,x)\,d\cS^\alpha(x)\,.
\end{equation}
\end{The}
In the sub-Riemannian framework, for distances with special symmetries and the proper choice of
the Riemannian surface meaure $\mu$, the density $\theta^\alpha(\mu,\cdot)$ is a geometric constant
that can be computed with a precise geometric interpretation. 
Then the previous formula is expected to have a potentially wide range of applications in the 
computation of the spherical Hausdorff measure of sets in the sub-Riemannian context.
Indeed, this project was one of the motivations for the present note.

Here we are mainly concerned with the purely metric area formula, therefore we limit ourselves 
to provide some examples of applications to the Heisenberg group, leaving details along with further 
developments for subsequent work.

Let $\Sigma$ be a $C^1$ smooth curve in $\H$ equipped with the sub-Riemannian distance $\rho$.
This distance is also called Carnot-Carath\'eodory distance, see \cite{Gr1} for the relevant definitions. 
Whenever a left invariant Riemannian metric $g$ is fixed on $\H$, we can associate $\Sigma$ with
its intrinsic measure $\mu_{SR}$, see \cite{Mag13Vit} for more details.
We will assume that $\Sigma$ has at least one {\em nonhorizontal point} $x\in\Sigma$, namely,
$T_x\Sigma$ is not contained in the horizontal subspace $H_x\H$, that is spanned by the
horizontal vector fields evaluated at $x$, \cite{Gr1}.

If we fix the size function $\zeta_{b,2}(S)=\diam(S)^2/4$ on closed balls and $x$ is nonhorizontal, 
then it is possible to compute explicitly $\theta^2(\mu_{SR},x)$, getting 
\[
\theta^2(\mu_{SR},x)=\alpha(\rho,g)\,,
\]
where $\alpha(\rho,g)$ is precisely the maximum among the lengths of all intersections of vertical lines
passing through the sub-Riemannian unit ball, centered at the origin. Here the length is computed with respect to
the scalar product given by the fixed Riemannian metric $g$ at the origin.
As an application of Theorem~\ref{the:metricspherical}, we obtain
\[
 \mu_{SR}=\alpha(\rho,g)\cS^2\res\Sigma\,,
\]
where $\cS^2$ is the spherical Hausdorff measure induced by $\zeta_{b,2}$.
The appearance of the geometric constant $\alpha(\rho,g)$ is a new phenomenon, due to the
use of Federer's density.
The nonconvex shape of the sub-Riemannian unit ball centered at the origin allows
$\alpha(\rho,g)$ to be strictly larger than the length $\beta(\rho,g)$ of the intersection 
of the same ball with the vertical line passing through the origin.
This feature of the sub-Riemannian unit ball shows that $\theta^2(\mu_{SR},x)$ and the the upper spherical density 
\[
\Theta^{*2}(\mu_{SR},x)=\limsup_{r\to0^+}\frac{\mu\ls \B(x,r)\rs}{r^2}
\]
differ. In fact, setting $t\in\ls \alpha(\rho,g),\beta(\rho,g)\rs$, we get
\[
\Theta^{*2}(\mu_{SR},x)=\beta(\rho,g)<t< \alpha(\rho,g)=\theta^2(\mu_{SR},x)\quad\mbox{for all}\quad x\in \cN,
\]
where $\cN=\lb x\in\Sigma: T_x\Sigma\ \mbox{is not horizontal} \rb$ and we also have
\begin{equation}\label{eq:alpha}
\mu_{SR}(\cN)=\alpha(\rho,g)\cS^2(\cN)>t\,\cS^2(\cN)\,.
\end{equation}
As a consequence of \eqref{eq:alpha}, in the inequality (1) of 2.10.19 in \cite{Federer69},
with $\mu=\mu_{SR}$ and $A=\cN$, the constant $2^m$ with $m=2$ cannot be replaced by one.
Moreover, even in the case we weaken the inequality (1) of 2.10.19 in \cite{Federer69} replacing 
the Hausdorff measure with the spherical Hausdorff measure, then \eqref{eq:alpha} still
shows that $2^m$ with $m=2$ cannot be replaced by one.
In the case $m=1$, it is possible to show, by an involved construction of a purely $(\cH^1,1)$
unrectifiable set of the Euclidean plane, that $2^m$ is even sharp,
see the example of 3.3.19 of \cite{Federer69}.
Somehow, our curve with nonhorizontal points has played the role of a more manageable unrectifiable set. 
Incidentally, the set $\cN$ is purely $(\cH^2,2)$ unrectifiable with respect to $\rho$, see \cite{AmbKir00}.

The connection between rectifiability and densities was already pointed out in \cite{PreTis92},
where the authors improve in a general metric space $X$ the upper estimate for $\sigma_1(X)$,
related to the so-called Besicovitch's $\frac{1}{2}$-problem.
According to \cite{PreTis92}, $\sigma_k(X)$ for some positive integer $k$ is the infimum among all 
positive numbers $t$ having the property that for each $E\subset X$ with $\cH^k(E)<+\infty$ 
and such that 
\[ 
\liminf_{r\to0^+}\frac{\cH^k\ls E\cap B(x,r)\rs}{c_k 2^kr^k}>t
\]
for $\cH^k$-a.e. $x\in E$ implies that $E$ is countably $k$-rectifiable, where
it is assumed that the open ball $B(x,r)$ has diameter equal to $2r$ for all $(x,r)\in X\times (0,+\infty)$ and
$\cH^k$ arises from the Carath\'eodory construction by the size function $\zeta(S)=c_k\,\diam(S)^k$.

If we equip $\H$ with the so-called Kor\'anyi distance $d$, see for instance Section~1.1 of \cite{KorRei95},
then a different application of Theorem~\ref{the:metricspherical} gives a lower estimate for $\sigma_2(\H,d)$.
In fact, we can choose $\Sigma_0$ to be a bounded open interval of the vertical line of $\H$ passing
through the origin. This set is purely $(\cH^2,2)$ unrectifiable. We define
$\zeta^d_{b,2}(S)=\diam_d(S)^2/4$ on closed balls, where the diameter $\diam_d(S)$ refers to the distance $d$,
and consider the intrinsic measure $\mu_{SR}$ of $\Sigma_0$.
By the convexity of the $d$-unit ball centered at the origin, the corresponding Federer density
$\theta^2_d(\mu_{SR},x)$ at a nonhorizontal point $x$ satisfies
\[
 \theta^2_d(\mu_{SR},x)=\alpha(d,g)\,,
\]
where $\alpha(d,g)$ is the length of the intersection of the Kor\'anyi unit ball cantered at the origin
with the vertical line passing through the origin. Following the previous notation, by Theorem~\ref{the:metricspherical}
we get 
\[
 \mu_{SR}=\alpha(d,g)\,\cS_d^2\res\Sigma_0\,,
\]
where $\cS^2_d$ is the spherical Hausdorff measure induced by $\zeta^d_{b,2}$. Since we have
\[
\lim_{r\to0^+}\frac{\cS_d^2\res\Sigma_0\ls B(x,r)\rs}{r^2}=
\lim_{r\to0^+}\frac{\mu_{SR}\ls B(x,r)\rs}{\alpha(d,g)r^2}=
\lim_{r\to0^+}\frac{\mu_{SR}\ls \B(x,r)\rs}{\alpha(d,g)r^2}
=1\,,
\]
and an easy observation shows that $\cS_d^2\res\Sigma_0\le 2\,\cH^2_2\res\Sigma_0$, 
it follows that 
\[
\frac{1}{2}=\lim_{r\to0}\frac{\cS_d^2\res\Sigma_0\ls B(x,r)\rs}{ 2r^2}
\le \liminf_{r\to0^+}
\frac{\cH^2_d\ls \Sigma_0\cap B(x,r)\rs}{r^2}\le\lim_{r\to0^+}
\frac{\cS^2_d\res\Sigma_0\ls B(x,r)\rs}{r^2}=1\,.
\]
This implies that $\sigma_2(\H,d)\ge 1/2$. 
Up to this point, we have seen how the geometry of the sub-Riemannian unit ball affects 
the geometric constants in estimates between measures. However, also the opposite direction is possible.
In fact, considering the previous subset $\cN$ and taking into account (1) of 2.10.19 in \cite{Federer69} with $m=2$, 
we get 
\[
\mu_{SR}(\cN)\le 4\, \beta(\rho,g)\,\cS^2(\cN)\,,
\]
hence the equality of \eqref{eq:alpha} leads us to the estimate
\[
1<\frac{\alpha(\rho,g)}{\beta(\rho,g)}\le 4\,.
\] 
It turns out to be rather striking that abstract differentiation theorems for measures 
can provide information on the geometric structure of the sub-Riemannian unit ball.
Precisely, we cannot find any left invariant sub-Riemannian distance $\tilde\rho$ in the Heisenberg group 
such that the geometric ratio $\alpha(\tilde\rho,g)/\beta(\tilde\rho,g)$ is greater than 4.
These facts clearly leave a number of related questions, so that the present note may
represent a starting point to establish deeper relationships 
between results of sub-Riemannian geometry and measure-theoretic results.

In particular, further motivations to study sub-Riemannian metric spaces may also arise 
from abstract questions of Geometric Measure Theory.
Clearly, to understand and carry out this demanding program more investigations are needed.


\begin{thebibliography}{99}



\bibitem{AmbKir00}{\sc L.~Ambrosio, B.~Kirchheim},
{\em Rectifiable sets in metric and Banach spaces}, Math. Ann. {\bf 318},
527-555, (2000)


\bibitem{Federer69}{\sc H.~Federer},
{\em Geometric Measure Theory}, Springer,  (1969)


\bibitem{Gr1}{\sc M.~Gromov},
{\em Carnot-Carath\'eodory spaces seen from within}, in {\em Subriemannian Geometry}, 
Progress in Mathematics, {\bf 144}. ed. by A. Bellaiche and J. Risler, Birkhauser Verlag, Basel, (1996).


\bibitem{Kir94} {\sc B.~Kirchheim},
{\em Rectifiable metric spaces: local
structure and regularity of the Hausdorff measure},
Proc. Amer. Math. Soc., {\bf 121}, 113-123, (1994).

\bibitem{KorRei95}{\sc A.~Kor\'anyi, H.~M.~Reimann},
{\em Foundation for the Theory of Quasiconformal Mappings
on the Heisenberg Group}, Adv. Math., {\bf 111}, 1-87, (1995).


\bibitem{Mag13Vit}{\sc V.~Magnani, D.~Vittone},
{\em An intrinsic measure for submanifolds in stratified groups},
J. Reine Angew. Math., {\bf 619}, 203-232, (2008) 


\bibitem{PreTis92}{\sc D.~Preiss, J.~Ti\v{s}er},
{\em On Besicovitch's $\frac{1}{2}$-problem}, J. London Math. Soc. (2) {\bf 45}, n.2, 279-287, (1992)


\bibitem{SawWhe92}{\sc E.~Sawyer, R.~L.~Wheeden}, 
{\em Weighted inequalities for fractional integrals on Euclidean and homogeneous spaces},
Amer. J. Math., {\bf 114}, n.4, 813-874, (1992)

\end{thebibliography}
\end{document}